\documentstyle[12pt]{amsart}

\newtheorem{theo}{\bf Theorem}[section]
\newtheorem{propo}[theo]{\bf Proposition}
\newtheorem{lemma}[theo]{\bf Lemma}

\newtheorem{coro}[theo]{\bf Corollary}

\newcommand{\inv}{{\rm inv}}
\newcommand{\cP}{{\cal P}_{nt}}

\begin{document}

\title[Expected number of inversions]{Expected number of inversions after a
sequence of random adjacent transpositions}
\author{Henrik Eriksson, Kimmo Eriksson and Jonas Sj{\"o}strand}
\address{NADA \\
   KTH \\
   SE-100 44 Stockholm, Sweden}
\email{henrik@@nada.kth.se}
\date{November 11, 1999}

\begin{abstract}
  In the evolution of a genome, the gene sequence is sometimes
rearranged, for example by transposition of two adjacent gene
blocks. In biocombinatorics, one tries to reconstruct these
rearrangement incidents from the resulting permutation. It seems
that the algorithms used are too effective and find a shorter
path than the real one. For the simplified case of adjacent 
transpositions, we give expressions for the expected number of
inversions after $t$ random moves. This average can be much
smaller than $t$, a fact that has largely been neglected so far.
\end{abstract}
\maketitle

\section{Background}

\noindent
The genome rearrangement problem is a combinatorial problem arising
in the area of molecular evolution. Basically, it can be stated as
a problem about permutations of gene sequences. 
Given a permutation $\pi$ as a word in the symbols $\{1,\dots,n\}$ 
(corresponding to genes), find the ``best'' path to the identity 
permutation when the feasible steps are block moves (removing a 
contiguous segment and inserting it somewhere else) and block
reversals 
(reversing the order of a segment). The shortest path between the two
permutations is the {\em parsimonious} solution, and finding
algorithms
for computing the shortest (or at least a short) path has been given a 
good deal of attention in \cite{Gu}, \cite{Er} etc.
Although no particular solution is more probable than the
parsimonious one, the shortest distance is not necessarily the most
probable length of a path. The number of shortest paths from the
identity to $\pi$ is much less than the number of paths using just a
few extra steps.  If a probabilistic model of the process is formulated,
a maximum-likelihood distance could be defined.
We have not seen this problem considered in the literature.

One approach is to try to determine the expected distance to the
identity permutation after a random walk of given length. This seems
to be a difficult problem. Of course, one can obtain intuition from
computer simulation, but a mathematical treatment would be
preferable.

In the present paper, we simplify the model so that the only steps 
allowed are adjacent transpositions. Observe that the set of adjacent
transpositions is the intersection of the set of block moves and the
set of block reversals. For this simpler problem, we are able to
obtain good lower and upper bounds by modelling the random walks in
terms of a discrete heat equation.

\section{Introduction}
\noindent
We are studying random walks of fixed length $t$ on the Cayley graph of
the Coxeter group $A_n$, which is isomorphic to the symmetric group
$S_{n+1}$.  The Cayley graph has an edge between two permutations
if one is obtained from the other by an adjacent transposition.
(See \cite{Al} for similar problems on $S_{n+1}$ and other groups.)

The walk starts at the identity permutation $1234\ldots$
and consists of $t$ random steps, chosen with uniform 
probability among the $n$ possible adjacent transpositions.  Let $\pi$ be
the permutation where this random walk stops and let $\inv(\pi)$ denote
the number of inversions. The shortest possible walk from the identity
to $\pi$ has length $\inv(\pi)$, so this number is less than or
equal to $t$. Clearly, for $t=1$, all permutations have $\inv(\pi)=1$,
but for $t\ge 2$, a later move may cancel an inversion 
created by an earlier move.
We would like to determine the expected
number of inversions $E(\inv(\pi))$, 
and we will denote it by $E_{nt}$ to make 
the dependence on the parameters $n$ and $t$ explicit. 

The set of adjacent transpositions is denoted by
$S=\{s_1,\dots,s_n\}$ where $s_i$ is the transposition of the
positions $i$ and $i+1$. Let $\cP$ be the set of all walks of length $t$, 
that is of all words in $S$ of length $t$: 
$$\cP = \{s_{i_1}s_{i_2}\dots s_{i_t} : 1\le i_1,\dots,i_t \le n\}.$$
Obviously $\cP$ has cardinality $n^t$. As the same notation is used for
a word and its product $\pi=s_{i_1}s_{i_2}\dots s_{i_t}$,  the notation
$\cP$ will be used also for the multiset of permutations.  
By counting all inversions in $\cP$ we can find the average number. 
$$n^t E_{nt} = \sum_{\pi\in\cP} \inv(\pi).$$

\noindent
Using a computer, we have calculated the integers $n^t E_{nt}$ 
for $n,t\le 10$. These data suggested a formula  of
unexpected simplicity.
Let $C_i$ denote the Catalan number $\frac{1}{i+1}{2i \choose i}$.

\begin{theo}\label{th:main}
For a fixed $t$ and for all $n\ge t$,
the expected number of inversions after $t$ random adjacent
transpositions is
$$E_{nt}= t-\frac{2}{n}{t\choose 2}+\sum_{r=2}^t \frac{(-1)^r}{n^r}
\left[2^r C_r {t\choose r+1} + 4d_r {t\choose r}\right] ,$$
where $d_2,d_3,d_4,\ldots$ is a certain integer sequence starting
by $0, 1, 9, 69, 510$. No expression for $d_r$ is known, but the 
following inequalities hold.

$$t-\frac{2}{n}{t\choose 2}+\sum_{r=2}^t \frac{(-1)^r}{n^r}
\left[2^r C_r {t\choose r+1} - 2^{r-1}C_{r-1} {t\choose r}\right]\le E_{nt}$$
$$E_{nt}\le t-\frac{2}{n}{t\choose 2}+\sum_{r=2}^t \frac{(-1)^r}{n^r}
2^r C_r {t\choose r+1}$$
\end{theo}
 
A direct proof of the theorem seems difficult, so
our approach has been a reformulation of the problem to a 
discrete heat flow model.

\section{The heat flow analogy}
\noindent
Instead of directly counting all inversions in $\cP$, we introduce
the following fine-grading.  Fixing $n$ and $t$, let
$$p_{ij} := {\rm Prob}(\pi_i < \pi_j)$$
for a random permutation $\pi\in\cP$.  Equivalently, 
$$n^t p_{ij}= \# \{\pi\in\cP : \pi_i < \pi_j\}.$$
Since every inversion is counted by one such class, we have
\begin{equation}\label{eq:Esum} 
E_{nt}=\sum_{i>j}p_{ij}.
\end{equation}
The matrices $(p_{ij})$ can be computed recursively.
For $t=0$, the set $\cP$ consists of the identity permutation
only, so $(p_{ij})$ has ones above the main diagonal and zeroes
below, as in the leftmost matrix of Figure 1.  The transformation
to $t=1$, $t=2$ and so on turns out to be a heat flow process.
The total heat is invariant, for $p_{ij}+p_{ji}=1$ by the law of
the Excluded Middle.  The main diagonal in the matrix 
may be left blank, as in the figure, or we may set all its entries to 
$1/2$ so that the rule $p_{ij}+p_{ji}=1$ is satisfied.

\begin{figure}[htb] 
\addtolength{\unitlength}{1.5pt}
\thicklines
\noindent
\begin{picture}(40,50)(10,10)
\small
\put(10,10){0}\put(13,11.2){\line(1,0){6}}
\put(20,10){0}\put(23,11.2){\line(1,0){6}}
\put(30,10){0}\put(33,11.2){\line(1,0){6}}
\put(40,10){0}\put(43,14){\line(1,1){6}}
\put(10,20){0}\put(13,21.2){\line(1,0){6}}
\put(20,20){0}\put(23,21.2){\line(1,0){6}}
\put(30,20){0}\put(33,24){\line(1,1){6}}
\put(50,20){1}
\put(10,30){0}\put(13,31.2){\line(1,0){6}}
\put(20,30){0}\put(23,34){\line(1,1){6}}
\put(40,30){1}\put(43,31.2){\line(1,0){6}}
\put(50,30){1}
\put(10,40){0}\put(13,44){\line(1,1){6}}
\put(30,40){1}\put(33,41.2){\line(1,0){6}}
\put(40,40){1}\put(43,41.2){\line(1,0){6}}
\put(50,40){1}
\put(20,50){1}\put(23,51.2){\line(1,0){6}}
\put(30,50){1}\put(33,51.2){\line(1,0){6}}
\put(40,50){1}\put(43,51.2){\line(1,0){6}}
\put(50,50){1}
\put(11.3,14){\line(0,1){5}}
\put(21.3,14){\line(0,1){5}}
\put(31.3,14){\line(0,1){5}}
\put(11.3,24){\line(0,1){5}}
\put(21.3,24){\line(0,1){5}}
\put(51.3,24){\line(0,1){5}}
\put(11.3,34){\line(0,1){5}}
\put(41.3,34){\line(0,1){5}}
\put(51.3,34){\line(0,1){5}}
\put(31.3,44){\line(0,1){5}}
\put(41.3,44){\line(0,1){5}}
\put(51.3,44){\line(0,1){5}}
\end{picture}
\begin{picture}(50,50)(0,10)
\small
\put(10,10){0}\put(13.6,11.2){\line(1,0){5}}
\put(20,10){0}\put(23.6,11.2){\line(1,0){5}}
\put(30,10){0}\put(33.6,11.2){\line(1,0){5}}
\put(40,10){$x$}\put(43.6,13){\line(1,1){6}}
\put(10,20){0}\put(13.6,21.2){\line(1,0){5}}
\put(20,20){0}\put(23.6,21.2){\line(1,0){5}}
\put(30,20){$x$}\put(33.6,23){\line(1,1){6}}
\put(48,20){1-$x$}
\put(10,30){0}\put(13.6,31.2){\line(1,0){5}}
\put(20,30){$x$}\put(23.6,33){\line(1,1){6}}
\put(38,30){1-$x$}\put(45,31.2){\line(1,0){4}}
\put(50,30){1}
\put(10,40){$x$}\put(13.6,43){\line(1,1){6}}
\put(28,40){1-$x$}\put(35,41.2){\line(1,0){4}}
\put(40,40){1}\put(43.6,41.2){\line(1,0){5}}
\put(50,40){1}
\put(18,50){1-$x$}\put(25,51.2){\line(1,0){4}}
\put(30,50){1}\put(33.6,51.6){\line(1,0){5}}
\put(40,50){1}\put(43.6,51.6){\line(1,0){5}}
\put(50,50){1}
\put(11.2,14.4){\line(0,1){4}}
\put(21.2,14.4){\line(0,1){4}}
\put(31.2,14.4){\line(0,1){4}}
\put(11.2,24.4){\line(0,1){4}}
\put(21.2,24.4){\line(0,1){4}}
\put(51.2,24.4){\line(0,1){4}}
\put(11.2,34.4){\line(0,1){4}}
\put(41.2,34.4){\line(0,1){4}}
\put(51.2,34.4){\line(0,1){4}}
\put(31.2,44.4){\line(0,1){4}}
\put(41.2,44.4){\line(0,1){4}}
\put(51.2,44.4){\line(0,1){4}}

\end{picture}
\begin{picture}(40,50)(0,10)
\scriptsize
\put(10.4,10){0}\put(13.6,11.0){\line(1,0){5}}
\put(20.4,10){0}\put(23.6,11.0){\line(1,0){5}}
\put(29.1,10){2$x^2$}\put(33.9,11.0){\line(1,0){2.6}}
\put(37.0,10){2$x${\tt -}3$x^2$}\put(43.2,13){\line(1,1){6}}
\put(10.4,20){0}\put(13.6,21.0){\line(1,0){5}}
\put(19.1,20){2$x^2$}\put(23.9,21.0){\line(1,0){2.6}}
\put(27.0,20){2$x${\tt -}4$x^2$}\put(33.2,23){\line(1,1){6}}
\put(43,20){1{\tt -}2$x${\tt +}3$x^2$}
\put( 9.1,30){2$x^2$}\put(13.9,31.0){\line(1,0){2.6}}
\put(17.0,30){2$x${\tt -}4$x^2$}\put(23.2,33){\line(1,1){6}}
\put(33,30){1{\tt -}2$x${\tt +}4$x^2$}\put(47,31.0){\line(1,0){2}}
\put(49.3,30){1{\tt -}2$x^2$}
\put( 7.0,40){2$x${\tt -}3$x^2$}\put(13.2,43){\line(1,1){6}}
\put(23,40){1{\tt -}2$x${\tt +}4$x^2$}\put(37,41.0){\line(1,0){2}}
\put(39.3,40){1{\tt -}2$x^2$}\put(47.3,41.0){\line(1,0){2.5}}
\put(50.5,40){1}
\put(13,50){1{\tt -}2$x${\tt +}3$x^2$}\put(27,51.0){\line(1,0){2}}
\put(29.3,50){1{\tt -}2$x^2$}\put(37.3,51.0){\line(1,0){2.5}}
\put(40.5,50){1}\put(43.6,51.0){\line(1,0){5}}
\put(50.5,50){1}
\put(11.3,14.4){\line(0,1){4}}
\put(21.3,14.4){\line(0,1){4}}
\put(31.3,14.4){\line(0,1){4}}
\put(11.3,24.4){\line(0,1){4}}
\put(21.3,24.4){\line(0,1){4}}
\put(51.3,24.4){\line(0,1){4}}
\put(11.3,34.4){\line(0,1){4}}
\put(41.3,34.4){\line(0,1){4}}
\put(51.3,34.4){\line(0,1){4}}
\put(31.3,44.4){\line(0,1){4}}
\put(41.3,44.4){\line(0,1){4}}
\put(51.3,44.4){\line(0,1){4}}

\end{picture}
\caption{Two heat flow steps}\label{fig:twosteps}
\end{figure}
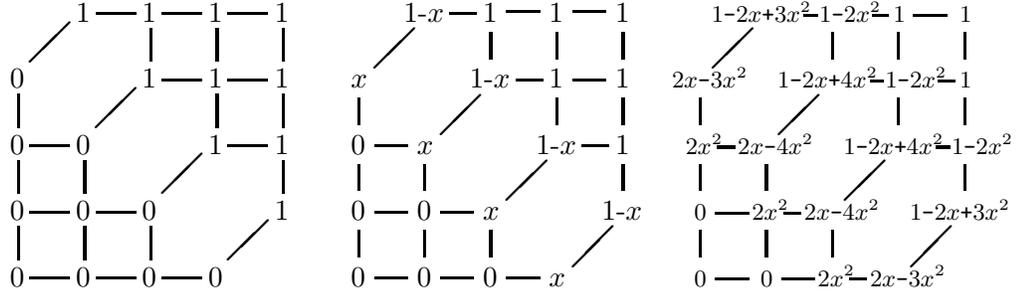

For any graph with real numbers (signifying temperature or heat) on the 
vertices, a 
{\em heat flow process} with thermal conductivity $x$ means the following.
In each step, every vertex sends the fraction $x$ of its heat
to each of its neighbours, at the same time receiving that
same fraction of its neighbours' heat.  Two steps are shown in
Fig.\,\ref{fig:twosteps}.

\begin{propo}\label{pr:heat flow} 
The sequence of $(p_{ij})$-matrices for $t=0,1,2,\dots$
describes a heat flow process with conductivity $x=\frac{1}{n}$ on the graph
depicted in Fig.1. The expected number of inversions, $E_{nt}$, equals the
total heat below the diagonal. 
\end{propo}
\begin{pf}
Consider the $(p_{ij})$-matrix after $t$ steps and let $\pi$ be
one of the permutations contributing to $p_{ij}$, that is,
$\pi$ satisfies $\pi_i < \pi_j$.  Each neighbour of $p_{ij}$
corresponds to a move that affects either $\pi_i$ or $\pi_j$.
For example, the neighbour to the right, $p_{i,j+1}$, corresponds 
to the transposition $(\pi_j,\pi_{j+1})$. When $i$ and $j$ are adjacent, 
the transposition $(\pi_i,\pi_j)$ is possible, which explains the graph 
edges across the main diagonal.

Except for these moves (at most four), the new $p'_{ij}$ would
be the same as the old $p_{ij}$, but now the following is true:
\begin{equation}\label{eq:neighbour}
p'_{ij}=p_{ij} + \frac{1}{n}\sum(p_{\rm neighbour}-p_{ij}),
\end{equation}
where the sum is taken over all graph neighbours of $p_{ij}$.
For after, say, the transposition $(\pi_j,\pi_{j+1})$, the 
$p'_{ij}$-condition $\pi'_i<\pi'_j$ means that we must have
had $\pi_i<\pi_{j+1}$, which is the $p_{i,j+1}$-condition.
\end{pf}

\noindent
From now on our analysis concerns the more general heat flow process
where $x$ is not necessarily $\frac{1}{n}$. For the matrix entries
we write $p_{ij}(x)$ and for the total heat below
the diagonal we use the notation $E_{nt}(x)$. For example,
Fig.\,\ref{fig:twosteps} demonstrates that
$E_{41}(x)=4x$ and $E_{42}(x)=8x-8x^2$. 

\begin{figure}[h]
\addtolength{\unitlength}{1.5pt}
\thicklines
\noindent
\begin{center}
\begin{picture}(60,50)(10,10)
\small
\put(10,10){0}\put(13,12){\line(1,0){6}}
\put(20,10){0}\put(23,12){\line(1,0){6}}
\put(30,10){0}\put(33,12){\line(1,0){6}}
\put(40,10){0}\put(43,12){\line(1,0){6}}
\put(50,10){$\frac{1}{2}$}
\put(10,20){0}\put(13,22){\line(1,0){6}}
\put(20,20){0}\put(23,22){\line(1,0){6}}
\put(30,20){0}\put(33,22){\line(1,0){6}}
\put(40,20){$\frac{1}{2}$}\put(43,22){\line(1,0){6}}
\put(50,20){1}
\put(10,30){0}\put(13,32){\line(1,0){6}}
\put(20,30){0}\put(23,32){\line(1,0){6}}
\put(30,30){$\frac{1}{2}$}\put(33,32){\line(1,0){6}}
\put(40,30){1}\put(43,32){\line(1,0){6}}
\put(50,30){1}
\put(10,40){0}\put(13,42){\line(1,0){6}}
\put(20,40){$\frac{1}{2}$}\put(23,42){\line(1,0){6}}
\put(30,40){1}\put(33,42){\line(1,0){6}}
\put(40,40){1}\put(43,42){\line(1,0){6}}
\put(50,40){1}
\put(10,50){$\frac{1}{2}$}\put(13,52){\line(1,0){6}}
\put(20,50){1}\put(23,52){\line(1,0){6}}
\put(30,50){1}\put(33,52){\line(1,0){6}}
\put(40,50){1}\put(43,52){\line(1,0){6}}
\put(50,50){1}
\put(11.3,15){\line(0,1){3}}
\put(21.3,15){\line(0,1){3}}
\put(31.3,15){\line(0,1){3}}
\put(41.3,15){\line(0,1){3}}
\put(51.3,15){\line(0,1){3}}
\put(11.3,25){\line(0,1){3}}
\put(21.3,25){\line(0,1){3}}
\put(31.3,25){\line(0,1){3}}
\put(41.3,25){\line(0,1){3}}
\put(51.3,25){\line(0,1){3}}
\put(11.3,35){\line(0,1){3}}
\put(21.3,35){\line(0,1){3}}
\put(31.3,35){\line(0,1){3}}
\put(41.3,35){\line(0,1){3}}
\put(51.3,35){\line(0,1){3}}
\put(11.3,45){\line(0,1){3}}
\put(21.3,45){\line(0,1){3}}
\put(31.3,45){\line(0,1){3}}
\put(41.3,45){\line(0,1){3}}
\put(51.3,45){\line(0,1){3}}
\end{picture}
\end{center}
\caption{Grid graph with initial values}\label{fig:grid}
\end{figure}
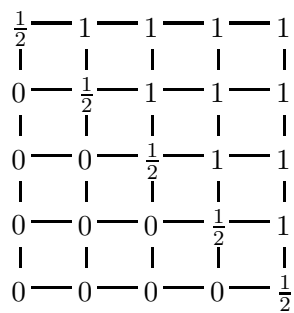

This analysis is complicated by the special edges across
the diagonal.
However, if we replace the graph of Fig.\,\ref{fig:twosteps} 
by the simple grid graph
of Fig.\,\ref{fig:grid} and set all diagonal values to $1/2$, then the heat flow
process is unchanged!  For thanks to the symmetry property
$p_{j,j+1}=1-p_{j+1,j}$ we have
$$p_{j,j+1}-p_{j+1,j}=1-2p_{j+1,j}= 2(1/2-p_{j+1,j}).$$
In other words, the loss of the neighbour across the diagonal is compensated
for by the two new neighbours on the diagonal.

\begin{propo}\label{pr:heat flow-grid} 
The sequence of $p_{ij}(x)$-matrices for $t=0,1,2,\dots$
describes a heat flow process on the 
$(n\!+\!1)\times (n\!+\!1)$ grid graph depicted in Fig.\,\ref{fig:grid}.
\end{propo}

\subsection{Hot boundary condition}
The above heat flow process on a grid with insulated boundary can be
reformulated as a heat flow process on the lower triangle
with the {\em hot boundary condition} $p_{ii}=1/2$ on the diagonal.
This should be obvious from the fact that the only connection
between the lower and upper triangle of the grid graph is the
diagonal, and the property of the original process that
the diagonal has constant temperature $p_{ii}=1/2$.

 Note that the subdiagonal
element $p_{j\!+\!1,j}(x)$ receives $2x/2$ from its diagonal neighbours
and sends back $2x\, p_{j\!+\!1,j}(x)$. The net heat transfer to the 
lower triangle
is $\sum_j (x\!-\!2x\,p_{j\!+\!1,j}(x))$, so we have the following result.
\begin{propo}\label{pr:Erecursion} 
$$ E_{n\,t\!+\!1}(x) =E_{n\,t}(x) + nx -2x\sum_j p_{j\!+\!1,j}$$
\end{propo}
\noindent
For example, in Fig.\,\ref{fig:twosteps} we see that
$ 8x-8x^2 = 4x + 4x -2x(4x).$

\subsection{Symmetric model and semi-infinite model}
The insulated left and lower boundaries can be gotten rid of in two
different ways, by symmetric extension or by just neglecting their effect.  

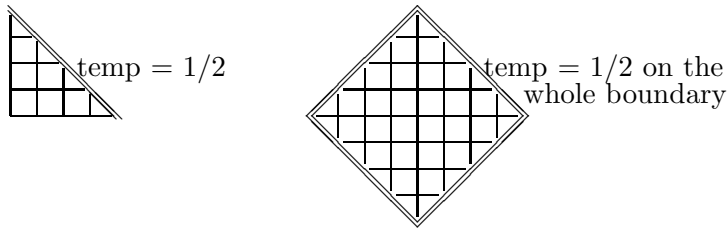
\begin{figure}[htb]
\begin{center}
\begin{picture}(110,40)(0,-30)
\put(0,0){\line(0,1){38}}
\put(10,0){\line(0,1){28}}
\put(20,0){\line(0,1){18}}
\put(30,0){\line(0,1){8}}
\put(0,0){\line(1,0){38}}
\put(0,10){\line(1,0){28}}
\put(0,20){\line(1,0){18}}
\put(0,30){\line(1,0){8}}
\put(-1,40){\line(1,-1){41}}
\put(-1,42){\line(1,-1){43}}
\small
\put(25,15){temp = 1/2}
\end{picture}
\begin{picture}(80,80)(-40,-30)
\put(0,-38){\line(0,1){76}}
\put(10,-28){\line(0,1){56}}
\put(20,-18){\line(0,1){36}}
\put(30,-8){\line(0,1){16}}
\put(-10,-28){\line(0,1){56}}
\put(-20,-18){\line(0,1){36}}
\put(-30,-8){\line(0,1){16}}
\put(-38,0){\line(1,0){76}}
\put(-28,10){\line(1,0){56}}
\put(-18,20){\line(1,0){36}}
\put(-8,30){\line(1,0){16}}
\put(-28,-10){\line(1,0){56}}
\put(-18,-20){\line(1,0){36}}
\put(-8,-30){\line(1,0){16}}
\put(0,40){\line(1,-1){40}}
\put(0,42){\line(1,-1){42}}
\put(0,40){\line(-1,-1){40}}
\put(0,42){\line(-1,-1){42}}
\put(0,-40){\line(1,1){40}}
\put(0,-42){\line(1,1){42}}
\put(0,-40){\line(-1,1){40}}
\put(0,-42){\line(-1,1){42}}
\small
\put(25,15){temp = 1/2 on the}
\put(40,5){whole boundary}
\end{picture}
\end{center}
\caption{Triangle with hot diagonal and symmetric extension of
the same problem}\label{fig:triangle}
\end{figure}

The reflection trick in Fig.\,\ref{fig:triangle}
demonstrates that the diamond graph with hot boundaries all around
is equivalent to the triangle. This model is what we will use for our exact
solution later in this paper. 

Neglecting the effect of the left and lower boundaries is equivalent
to moving them to infinity. Then we are left with the whole half-plane
below an infinite hot diagonal. As we will see, this problem is not so
difficult and we can use its solution as a lower bound.  

It is clear that the temperature at a given inner point in the
diamond model must be at least as hot as a point at the same distance
from the diagonal in the semi-infinite
model, since the former point has heat flowing to it from three additional
sides. Hence, solving the semi-infinite model gives a lower bound for 
the actual finite case.

\section{Heat flow combinatorics}
\noindent
Our goal in this section is to find combinatorial expressions
for the $p_{ij}(x)$-matrices describing heat flow on a 
triangular grid graph. As explained above, the triangular graph can be
considered as embedded either in a finite diamond graph or in
a semi-infinite grid graph.  

At day $t=0$ the diagonal entries are $p_{ii}=\frac{1}{2}$ with zeroes
below the diagonal.
Referring to Fig.\,\ref{fig:twosteps}
and the recursion ($p'$ means next day)
\begin{equation}
p'_{ij}(x)=p_{ij}(x) + x\sum\left( p_{\rm neighbour}(x)-p_{ij}(x)\right),
\end{equation}
 it is obvious
that the entries $p_{ij}(x)$ in step $t$ will be polynomials in
$x$ of degree $t$ (or less).  We can give each coefficient in
these polynomials a combinatorial significance. Vaguely expressed,
they count journeys for $t$ days from the hot boundary to the
location of $p_{ij}(x)$. The recursion states that such a journey 
ending in a certain location a certain day [$p'_{ij}(x)$] 
may have been either at the same location yesterday
[$p_{ij}(x)$] and had a resting day or at a neighbouring location
[$p_{\rm neighbour}(x)$] and had a travel day (the $x$-factor means
travel day) or at the same location [$p_{ij}(x)$] and travelled
half-way to the neighbour and then back (the $-x$-factor means round trip).

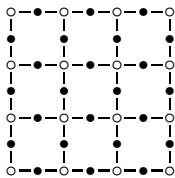
\begin{figure}[htb]
\begin{center}
\begin{picture}(60,50)(10,10)
\put(0,0){\circle{3}}
\put(10,0){\circle*{3}}
\put(20,0){\circle{3}}
\put(30,0){\circle*{3}}
\put(40,0){\circle{3}}
\put(50,0){\circle*{3}}
\put(60,0){\circle{3}}
\put(0,20){\circle{3}}
\put(10,20){\circle*{3}}
\put(20,20){\circle{3}}
\put(30,20){\circle*{3}}
\put(40,20){\circle{3}}
\put(50,20){\circle*{3}}
\put(60,20){\circle{3}}
\put(0,40){\circle{3}}
\put(10,40){\circle*{3}}
\put(20,40){\circle{3}}
\put(30,40){\circle*{3}}
\put(40,40){\circle{3}}
\put(50,40){\circle*{3}}
\put(60,40){\circle{3}}
\put(0,60){\circle{3}}
\put(10,60){\circle*{3}}
\put(20,60){\circle{3}}
\put(30,60){\circle*{3}}
\put(40,60){\circle{3}}
\put(50,60){\circle*{3}}
\put(60,60){\circle{3}}
\put(0,50){\circle*{3}}
\put(20,50){\circle*{3}}
\put(40,50){\circle*{3}}
\put(60,50){\circle*{3}}
\put(0,30){\circle*{3}}
\put(20,30){\circle*{3}}
\put(40,30){\circle*{3}}
\put(60,30){\circle*{3}}
\put(0,10){\circle*{3}}
\put(20,10){\circle*{3}}
\put(40,10){\circle*{3}}
\put(60,10){\circle*{3}}
\put(0,3){\line(0,1){4}}
\put(0,13){\line(0,1){4}}
\put(0,23){\line(0,1){4}}
\put(0,33){\line(0,1){4}}
\put(0,43){\line(0,1){4}}
\put(0,53){\line(0,1){4}}
\put(20,3){\line(0,1){4}}
\put(20,13){\line(0,1){4}}
\put(20,23){\line(0,1){4}}
\put(20,33){\line(0,1){4}}
\put(20,43){\line(0,1){4}}
\put(20,53){\line(0,1){4}}
\put(40,3){\line(0,1){4}}
\put(40,13){\line(0,1){4}}
\put(40,23){\line(0,1){4}}
\put(40,33){\line(0,1){4}}
\put(40,43){\line(0,1){4}}
\put(40,53){\line(0,1){4}}
\put(60,3){\line(0,1){4}}
\put(60,13){\line(0,1){4}}
\put(60,23){\line(0,1){4}}
\put(60,33){\line(0,1){4}}
\put(60,43){\line(0,1){4}}
\put(60,53){\line(0,1){4}}
\put(3,0){\line(1,0){4}}
\put(13,0){\line(1,0){4}}
\put(23,0){\line(1,0){4}}
\put(33,0){\line(1,0){4}}
\put(43,0){\line(1,0){4}}
\put(53,0){\line(1,0){4}}
\put(3,20){\line(1,0){4}}
\put(13,20){\line(1,0){4}}
\put(23,20){\line(1,0){4}}
\put(33,20){\line(1,0){4}}
\put(43,20){\line(1,0){4}}
\put(53,20){\line(1,0){4}}
\put(3,40){\line(1,0){4}}
\put(13,40){\line(1,0){4}}
\put(23,40){\line(1,0){4}}
\put(33,40){\line(1,0){4}}
\put(43,40){\line(1,0){4}}
\put(53,40){\line(1,0){4}}
\put(3,60){\line(1,0){4}}
\put(13,60){\line(1,0){4}}
\put(23,60){\line(1,0){4}}
\put(33,60){\line(1,0){4}}
\put(43,60){\line(1,0){4}}
\put(53,60){\line(1,0){4}}
\end{picture}
\end{center}
\caption{Grid graph with mid-vertices}\label{fig:midgrid}
\end{figure}

\noindent
In order to make our statements precise, we will modify the grid
graph as in Fig.\,\ref{fig:midgrid}.  Each of the
original edges is split in two by a new {\em mid-vertex}. Mid-vertices
are introduced for counting purposes only, they do not carry heat.

Each day, the vertices on the hot boundary send out heat packets
with the value $\frac{x}{2}$ to their neighbours. These packets
are sent on and on, back and forth, always multiplied by $x$ 
or $-x$. Consider one such heat packet at a certain location in the
morning of a certain day. What can happen to it during the day? 
\begin{enumerate}
\item It stays on its vertex unchanged.
\item It travels a half-edge, gets multiplied by $x$, and travels
the other half-edge to the next vertex.
\item It travels a half-edge, gets multiplied by $-x$, and returns
the same half-edge to the same vertex.
\end{enumerate}
If the start value is $\frac{1}{2}$ and the journey has $r$ travel 
days (type 2 or 3), the final value is $\pm \frac{1}{2}x^r$. 
The sign depends on the number 
of days of type 3, and is easily seen to be $(-1)^{r+i-j}$
if the journey ends at 
$(i,j)$. Hence we have the following result.
\begin{lemma}\label{lm:coefficient}
The coefficient of the
$x^r$-term in $p_{ij}(x)$ is $(-1)^{r+i-j}$ times the
number of journeys from the hot boundary to $(i,j)$ in $t$ days, $r$ of
which are travel days.
\end{lemma}

\subsection{The semi-infinite model}\label{sec:inf}
Our next step is counting the journeys specified in the lemma. This is easy
in the semi-infinite case where all points on the $k$-subdiagonal are 
equivalent. We define the {\it sublevel} of $(i,j)$ as $k=i-j$. The journey
from sublevel $k$ to sublevel $0$ can be specified by three items:
\begin{itemize}
\item Out of the $t$ days, $r$ travel days must be chosen. This can be
done in $t \choose r$ ways.
\item For each of the $r$ travel days, horizontal or vertical travel must
be chosen. This can be done in $2^r$ ways. 
\item A Catalan walk in $2r$ half-steps from sublevel $k$ to sublevel 0
must be specified.
By Catalan walk we mean that sublevel 0 must not be reached until the last
half-step. It is well-known that the number of such walks is 
$$ {2r\!-\!1\choose r\!-\!k}-{2r\!-\!1\choose r\!-\!k\!-1}.$$
For $k=1$ this is the Catalan number $C_r$.
\end{itemize}

\noindent
Combining this journey count with Lemma \ref{lm:coefficient}, we
get the expressions for $p_{ij}(x)$ and eventually the total heat in
the triangle.
\begin{propo}
Ater $t$ time steps in the semi-infinite model,
\begin{equation}\label{eq:pinf}
p_{i,j}(x)={1\over 2}\sum_{r=k}^t (-1)^{r+k}{t \choose r}2^r
\left[ {2r\!-\!1\choose r\!-\!k}-{2r\!-\!1\choose r\!-\!k\!-1}\right] x^r,
\end{equation}
where $k=i-j$. The total heat under the diagonal is 
\begin{equation}\label{eq:Einf}
E_{nt}(x)=ntx-\frac{n+1}{2}\sum_{r=2}^t (-1)^{r}{t \choose r}2^r C_{r\!-\!1}x^r
\end{equation}
\end{propo}
\begin{pf}
  As there are $n\!+\!1\!-\!k$ locations on sublevel $k$, we must do the following sum.
$$ E_{nt}(x)=\sum_{k=1}^n (n\!+\!1\!-\!k)\,{1\over 2}\sum_{r=k}^t (-1)^{r+k}{t \choose r}2^r
\left[ {2r\!-\!1\choose r\!-\!k}-{2r\!-\!1\choose r\!-\!k\!-1}\right] x^r$$.
$$=\sum_{r=1}^t (-1)^{r}{t \choose r}2^{r} x^r\,
\sum_{k=1}^r (-1)^{r}\frac{n\!+\!1\!-\!k}{2} 
\left[ {2r\!-\!1\choose r\!-\!k}-{2r\!-\!1\choose r\!-\!k\!-1}\right],$$
and the last sum simplifies to give the desired result.
\end{pf}

\noindent
As we observed in Sec. \ref{sec:inf}, the total heat in the
semi-infinite model gives a lower bound for the total heat in
the finite case. 
In particular, we can plug in $x=\frac{1}{n}$ to obtain a lower bound for 
the expected number of inversions.
\begin{coro}
  The lower bound for $E_{nt}$ in Theorem \ref{th:main} holds true.
\end{coro}
\begin{pf}
  Substitute $x=\frac{1}{n}$ and collect like powers of $\frac{1}{n}$. 
\end{pf}

\noindent
Remarkably, the semi-infinite model also provides an upper bound for
the total heat in the finite case. By iterating the recursion 
\ref{pr:Erecursion} for the $E_{nt}$ in the finite case,
we obtain the following formula. 
\begin{lemma}
  Let $e_t(x)$ denote the sum of the subdiagonal entries of
  the $p_{ij}(x)$-matrix for time step $t$ in the finite case.  Then
  $$E_{nt}=ntx-{2x}[e_{t-1}(x)+e_{t-2}(x)+\cdots+e_{1}(x)].$$
\end{lemma}
\noindent
We know that all $p_{ij}(x)$ in the semi-infinite model are less than 
or equal to the $p_{ij}(x)$ in the finite case. In particular, the
subdiagonal sums must be less, so if we use them for $e_t(x)$ in the
lemma above, we obtain an upper bound for $E_{nt}$.
\begin{coro}
  The upper bound for $E_{nt}$ in Theorem \ref{th:main} holds true.
\end{coro}
\begin{pf}
  Use the lemma together with Eq. \ref{eq:pinf}, then substitute 
$x=\frac{1}{n}$ and simplify.
\end{pf}

\subsection{The finite case}
Lemma \ref{lm:coefficient} tells us that $p_{ij}(x)$ is
an alternating polynomial and that its coefficients counts
journeys from the hot boundary to $(i,j)$. In the finite case,
counting journeys is difficult when $n<t$, for there are four
boundaries and wherever you start it is possible to reach more than      
one of them in $t$ days. But when $n\ge t$, the situation is
better. 

The expression for the number of journeys with $r$ travel days
starting at $(i,j)$ and ending at a hot boundary used to be
$${t \choose r}2^r
\left[ {2r\!-\!1\choose r\!-\!k}-{2r\!-\!1\choose r\!-\!k\!-1}\right]$$
but for some $(i,j)$ that are close to two boundaries, this number will
now increase. For the extra journeys, horizontal and vertical steps
cannot be chosen freely, so the factor $2^r$ does not apply. For example,
from $(2,1)$ it is possible to reach the left boundary in two steps, at
least one of which must be horizontal. Therefore, the contribution 
of the extra journeys to $E_{nt}(x)$ will be of the form
\begin{equation}\label{eq:correction}
\sum_{r=2}^t (-1)^{r}{t \choose r}2g_r x^r.
\end{equation}
The factor 2 comes from symmetry. More important than the exact
value of the $g_r$-numbers is the fact that they do not depend on $n$.
Therefore, these correction terms become less and less important as
$n$ increases.  
For the expected number of inversions, $E_{nt}$, the correction terms
may be written as in Theorem \ref{th:main}.

\section{Open problems}
\begin{enumerate}
\item Is there a nice expression for the $d_r$-numbers of Theorem 
\ref{th:main}? 
\item Is there a nice expression for the $g_r$-numbers of Eq.
\ref{eq:correction}?
\item Can the analysis be extended to adjacent block tranpositions?
\item Can the analysis be extended to block reversals?
\item If the result of some random moves is a permutation with a
certain number of inversions, what number of moves is the most probable? 
\end{enumerate}

\end{document}